\magnification=1200

\centerline{Uniqueness in Rough Almost Complex Structures,}

\centerline{and Differential Inequalities.}
\bigskip\centerline{by Jean-Pierre Rosay.}
\bigskip\bigskip\noindent
INTRODUCTION. \hfill\break
The open unit disc in ${\bf C}$ will be denoted by ${\bf D}$.
Recall that an almost complex structure on ${\bf C}^n$ consists in
having for each $p\in {\bf C}^n$ an endomorphism $J=J(p)$ of the (real) tangent
space to ${\bf C}^n$ at $p$ satisfying $J^2=-{\bf 1}$, and that a map $u:{\bf D}\to
({\bf C}^n,J)$ is $J$-holomorphic if ${\partial u\over\partial y}(z)= [J(u(z))]\big({\partial u\over \partial x}(z)\big)$.\medskip\noindent
The main results in this paper are the following ones:
\medskip\noindent
{\bf Proposition 1.} {\it Let $J$ be a H\"older continuous
${\cal C}^{1\over 2}$ almost complex structure
defined on ${\bf C}^n$. Let $u:{\bf D}\to {\bf C}^n$ be a $J$-holomorphic map
(so $u$ is of class ${\cal C}^{1,{1\over 2}}$). If $u=0$ on some non empty open
subset of ${\bf D}$, then $u\equiv 0$.}
\bigskip\noindent
{\bf Proposition 2.} {\it There exists a smooth map $u:{\bf D}\to {\bf C}^2$ satisfying
$|{\partial u\over\partial \overline z}|\leq \epsilon (z)|{\partial u\over \partial z}|$
with $\epsilon (z)\to 0$ as $z\to 0$, and such that:

$u$ vanishes to infinite order at 0,

$u$ has a non-isolated zero at 0,

\noindent
but $u$ is not identically 0 near 0.}
\bigskip\noindent
Note that the phenomena of Proposition 2 cannot occur for scalar valued maps $u$ 
(i.e. for $u:{\bf D}\to {\bf C}$), see the Appendix. So, here we see a
difference between vector valued maps and functions.
\bigskip\noindent
I shall now explain the motivations for these questions and how they are related.
\hfill\break
It is well known and extremely easy to see that for H\"older continuous almost
complex structures, that are not Lipschitz, there is not an equivalent of unique analytic
continuation. Two distinct $J$ holomorphic maps $u$ and $v$ can agree on 
a non empty open subset of ${\bf D}$. Indeed, set $u(z)=(z,0)$, and $v(z)=(z,0)$
if ${\rm Im}~z\leq 0$, but $v(z)=(z,({\rm Im}~z)^k)$ if ${\rm Im}~z>0$. For $u$
and $v$ to be both $J$-holomorphic maps we simply need that $[J(z,0)](1,0)=(i,0)$,
and $[J(x+iy,y^k)](1,0)=(i,ky^{k-1}))$, if $y>0$. It is immediate that such $J$'s of
H\"older class ${\cal C}^{k-1\over k}$ can be defined.
This failure of uniqueness is essentially linked to the failure
of uniqueness for O.D.E. such as $y'=|y|^\alpha$, $\alpha <1$, for which the Lipschitz
condition is not satisfied. The possibly surprising fact is that uniqueness
holds when one of the maps is constant, at least for almost complex structures of class
${\cal C}^{1\over 2}$. 
\bigskip\noindent
If $J$ is close enough to the standard complex structure $J_{st}$, $J$-holomorphic maps
are characterized by an equation:
$${\partial u\over \partial \overline z}=
\overline{Q(u){\partial u\over \partial z}}~~~~~~(E)$$
where for $p\in {\bf C}^n$, $Q(p)$ is a ${\bf C}$-linear map, $Q(p)=0$ if $J(p)=J_{st}$,
and $Q$ has the same H\"older or ${\cal C}^k$ regularity as $J$.
\bigskip\noindent
Equation (E) leads to inequalities, in particular:
\bigskip
(1) If $Q(0)=0$ and $Q$ is Lipschitz continuous, (locally at least) one gets:
$$|{\partial u\over \partial \overline z}|\leq C|u|~~~~~(IN1)$$
where by rescaling the size of $C$ is irrelevant.
\bigskip
(2) If the operator norm of $Q(.)$ is $\leq \epsilon$,
$$|{\partial u\over\partial \overline z}|\leq \epsilon |{\partial u\over \partial z}|~~~~~~(IN2)$$
The size of $\epsilon$ is important.
\bigskip\noindent
In almost complex analysis, it is interesting to know which properties follow
from the equations and which one are merely consequences of the above inequalities.
(IN1) is very easy to use and there is a summary in the Appendix. (IN2) is useful
to get energy estimates (see Lemma 2.4.2 in [14], 
Remark 2 in 1.d in [8], 2.2 in [4]). Another differential inequality,
on the Laplacian, is used in [12] page 44, also to get energy estimates.
Contrary to (IN1), (IN2) apparently does
not yield good uniqueness results, as Proposition 2 illustrates. A much strengthened
version of (IN2) will be used for proving Proposition 1.
\bigskip\noindent
{\bf Proposition 3.} {\it Let $u:{\bf D}\to {\bf C}^n$ be a ${\cal C}^1$ map such that
for some $K>0$
$$|{\partial u\over \partial\overline z}|\leq K|u|^{1\over 2}|{\partial u\over \partial z}|~.$$
If $u$ vanishes on some non empty open subset of ${\bf D}$, then $u\equiv 0$.} 
\bigskip\noindent
{\bf Comments.}\hfill\break
1) The problems of uniqueness

(a) vanishing on an open set implies vanishing 

(b) vanishing to infinite order at a point implies vanishing
on a neighborhood \hfill\break
are not strictly related since we are not dealing with smooth maps.
\smallskip\noindent
2) The results in this paper are only partial results. More questions are raised than
solved, e.g. :
\smallskip Can Proposition 1 be extended to all H\"older continuous almost
complex structures?
\smallskip
In Proposition 1, can one replace vanishing on an open set
by vanishing to infinite order at some point, or having a non-isolated zero? 
(For non-H\"olderian structures, see an example in Remark 3).
\smallskip
And vice versa for Proposition 2.
\bigskip\noindent
The gap is huge between the positive result of Proposition 3 and
the counterexample of Proposition 2, in which (with the notation of
Proposition 2) we will have $|\epsilon (2^{-n})| = O({1\over n})$.
\bigskip\noindent
3) The proof of Proposition 1 will follow very closely the first, easier, steps
of the proof of Theorem 17.2.1 in [6], and will conclude with the arguments in
8.5 and 8.6 in [5],  mentioned in [6] page 10. A special case of this Theorem, 
gives the following: \hfill\break
Let $g$ be a ${\cal C}^1$ function defined on ${\bf D}$. Assume that
${\partial u\over \partial \overline z}+a(z){\partial g\over \partial z}=0$,
where $a$ is a ${\cal C}^1$, with $|a|<1$ . If $g$ vanishes on some non empty 
open subset of ${\bf D}$, then $g\equiv 0$.
\smallskip\noindent
For getting uniqueness for $J$-holomorphic discs there are very serious difficulties
preventing us from adapting the result:

- The proof makes crucial use of the differentiability of $a$ that we certainly do not have.

- We would need a generalization to vector valued maps $g$, with $a(z)$ operator valued.
The proof in [6] seems to be difficult to adapt if one does not assume $a(z)$ to be
a normal operator ($a^\ast a=aa^\ast$). Is such a generalization true?

- And first, and possibly worst: equation (E) is not a ${\bf C}$-linear equation.
\bigskip\noindent
4) Aronszajn type Theorems, such as Theorem 17.2.6 in [6], that prove that vanishing 
to infinite order implies vanishing, are difficult and the situation is very delicate,
as shown both by the fact that these results do not generalize to equations order $>2$ 
and by the second order counterexamples provided by Alinhac [1]. However, in case of
the differential inequality $|{\partial u\over \partial \overline z}|\leq C|u|$ (vector valued
case), the proof of uniqueness under the hypothesis of vanishing to infinite order
by using Carleman weights simplifies enormously, since one can take advantage of the
commutation of ${\partial\over \partial\overline z}$ and multiplication by
${1\over z^N}$, for introducing the weights ${1\over |z|^N}$, at `no cost'
(proving an estimate $\int_D|{\partial u\over \partial\overline z}|^2{1\over |z|^{2N}}
\geq K\int_D|u|^2{1\over |z|^{2N}}$, for $u\in {\cal C}^1_0(D)$, with $K>0$ independent
of $N$), avoiding thus
all the difficult part of the proof of Aronszajn's Theorem.
\vskip.5truein
\noindent
\S 1. Standard $\overline \partial$ estimates.
\smallskip\noindent
Notation: For the whole paper, we set $\phi (z)=\phi (x) =x+{x^2\over 2}$,
$z=x+iy$. So $\phi (z)<0$ if $-1\leq x<0$, $\phi$ is increasing on
$[-1,1]$, and $\phi''=1$.
\smallskip\noindent Integration on ${\bf D}$ will be simply denoted by
$\int$, standing for $\int_{\bf D}$.
\bigskip\noindent
The estimate of Lemma 1 is completely standard and its proof is included
only for the convenience of the reader. It is more usual to see the 
dual estimate used for solving $\overline \partial$ requiring the opposite sign
of the exponent in the weight function. For generalizations see [6], proof
of Theorem 17.2.1. For an elementary introduction, see 4.2 in [7] 
4.2 in Chapter IV. 
\bigskip\noindent
1.1. \hfill\break
{\bf Lemma 1.} {\it For any ${\cal C}^1$ function (or ${\bf C}^n$-valued map)
$v$ with compact support in ${\bf D}$,
and $\tau\geq 1$
$$\int |{\partial v\over \partial \overline z}|^2e^{\tau \phi}\geq
{1\over 10\tau} \int |{\partial v\over \partial z}|^2e^{\tau \phi}
+{\tau\over 20}\int|v|^2e^{\tau \phi}~.$$}
\bigskip\noindent
Proof. For short set $\overline \partial ={\partial \over \partial\overline z}$
and $\partial = {\partial \over \partial z}$. Consider the Hilbert space of measurable
functions $f$ on ${\bf D}$ such that $\int |f|^2e^{\tau\phi} <+\infty$, with the scalar product
$<f,g>_\phi=\int f\overline ge^{\tau \phi}$.
\smallskip\noindent
Let $\overline \partial^\ast$ be the adjoint of $\overline\partial$ in that space.
Elementary computations show that 
$$\overline \partial^\ast v = -\big( \partial v + {1\over 2} \tau \phi'\big)~.$$ 
One has
$$\int|\overline \partial v|^2e^{\tau \phi}=<\overline\partial v,\overline\partial v>_\phi =
<v,\overline \partial^\ast \overline\partial v>_\phi =
<v, \overline \partial~ \overline \partial^\ast v>_\phi +
<v, [\overline \partial^\ast ,
\overline \partial ]v>_\phi
=$$
$$
<\overline\partial^\ast v,  \overline \partial^\ast v>_\phi +
<v, [\overline \partial^\ast ,
\overline \partial ]v>_\phi
~.$$ 
Another immediate computation gives $[\overline\partial^\ast, \overline\partial] v=
{1\over 4}\phi''v={1\over 4}v$. (This is where positivity comes, from the convexity of
$\phi$).
Therefore $$\int |\overline \partial v|^2e^{\tau\phi}=
\int |\overline \partial^\ast v|^2e^{\tau\phi}+{\tau\over 4} 
\int |v|^2e^{\tau \phi}~\geq
{1\over 5\tau}\int |\overline\partial^\ast v|^2e^{\tau\phi}+
{\tau\over 4}\int |v|^2e^{\tau \phi}~.$$
\medskip\noindent
To finish the proof, it is enough to use, in the above inequality, the estimate from below
(simply using $(a+b)^2\geq {a^2\over 2}-b^2$):\hfill\break
$|\overline\partial^\ast v|^2=|\partial v+{1\over 2}\tau\phi' v|^2\geq
{1\over 2} |\partial v|^2 - {1\over 4}|\tau \phi' v|^2\geq
{1\over 2} |\partial v|^2- \tau^2 |v|^2$, since $0\leq \phi'\leq  2$.

\bigskip

\noindent
1.2. The next Lemma is narrowly tailored to the application in view. 
The hypotheses are made just ad-hoc and
we state the conclusion just as we shall need it, dropping an ${1\over \tau}\int|\nabla u|^2$
term on the right hand side, that we could keep.
\bigskip\noindent
{\bf Lemma 2.} {\it Let $v$ be a ${\cal C}^1$ map from ${\bf D}$ to ${\bf C}^n$, with
$|v|\leq 1$. If $w$ is a measurable function on ${\bf D}$ satisfying 
$|w(z)|\leq \theta |v|^{1\over 2}~{\rm Min}~(|{\partial v\over \partial z}|,1)$,
with $\theta={1\over 10}$.
then, for any $\tau \geq 1$:
$$\int|{\partial v\over \partial \overline z}-w(z)|^2e^{\tau \phi}
\geq {\tau\over 80}\int|u|^2e^{\tau \phi}~.$$}
\bigskip\noindent
Proof. We have
$$\int|{\partial v\over \partial \overline  z}-w(z)|^2e^{\tau\phi}\geq
{1\over 2}\int |{\partial v\over\partial\overline z}|^2e^{\tau\phi}-\int |w|^2e^{\tau\phi}$$
$$\geq {1\over 20\tau}\int|{\partial v\over \partial z}|^2e^{\tau\phi}
+{\tau\over 40}\int |v|^2e^{\tau\phi}
-\theta^2\int|v| [{\rm Min}~(|{\partial v\over \partial z}|,1)]^2e^{\tau\phi}~.$$
The Lemma will be established if we have the (much better than needed) point-wise estimate
$$\theta^2|v| [{\rm Min}~(|{\partial v\over \partial z}|,1)]^2
\leq {1\over 20\tau}|{\partial v\over \partial z}|^2+{\tau\over 80}|v|^2~.$$
At points where $|v|\leq {1\over \tau}$, the estimate is trivial.  
At points where $|v|\geq {1\over \tau}$, we simply use 
$\theta^2|v| [{\rm Min}~(|{\partial v\over \partial z}|,1)]^2\leq \theta^2 |v|$.
We have $\theta^2 |v|= \theta^2 |v|^{-1}|v|^2\leq \theta^2\tau |v|^2\
\leq {1\over 80}\tau |v|^2$. Q.E.D.
\bigskip\noindent
{\bf Remark 1.} The proof of Lemma 2 ends with a simple pointwise estimate,
because we have no regularity assumptions on $w$ (so, no helpful integration by parts
seems to be possible). The H\"older exponent 
${1\over 2}$ in Proposition 1,
will lead below to the consideration of a perturbation term
$w(z)$ satisfying the hypotheses of Lemma 2. That exponent ${1\over 2}$
seems to be the limit of what our approach can reach.
In Lemma 1, a big loss was taken. At some point, we wrote
$\int |\overline\partial^\ast v|^2e^{\tau\phi}
\geq {1\over 5\tau}\int |\overline\partial^\ast v|^2e^{\tau\phi}$.
This is just to say that we could add
${1\over 2} \int |\overline\partial^\ast v|^2e^{\tau\phi}$,
to the right hand side of the inequality in Lemma 1. However, as we shall see,
this would not solve the difficulty that we now explain.
If we now consider an almost complex structure of class ${\cal C}^\alpha$,
we are led to consider a perturbation term $w$ of the size
of $|v|^\alpha|{\partial v\over \partial z}|$. If, for simplicity, we assume
$|{\partial v\over \partial z}|\leq 1$, for finishing the proof of 
Lemma 2, under the hypothesis $|w|\leq \theta |v|^\alpha|{\partial v\over \partial z}|$, 
we would need an inequality of the type
$$|{\partial v\over \partial z}+\tau \phi'v|^2+{1\over \tau}|{\partial v\over
\partial z}|^2+\tau |v|^2\geq C |v|^{2\alpha}|{\partial v\over \partial z}|^2~,$$
where $C>0$ should be independent of $\tau$, and the first term on the left hand side
is to try to get advantage from the term previously dropped (that makes the 
second term superfluous, as it has been seen the proof of Lemma 1).
Suppose that at some point $z$, with $|z|<1$ and ${\rm Re}~z\geq -{1\over 2}$
(so ${1\over 2}\leq \phi'\leq 2$), we have:
${\partial v\over \partial z}+\tau \phi'v=0$, and $|v(z)|={\epsilon\over \tau}$,
with $0<\epsilon<{1\over 2}$
(so $|{\partial v\over\partial z}(z)|\leq 2\epsilon<1$). Then, the left hand side,
in the above inequality, is
$\tau (1+\phi'^2)|v|^2\leq  5\tau |v|^2= {5\epsilon^2\over \tau}$, while the right hand side is
$\geq {C \tau^2\over 4} |v|^{2+2\alpha}={C\epsilon^{2+2\alpha}\over 4}{1\over \tau^{2\alpha}}$.
For $\tau$ large, the 
inequality will not hold if $\alpha < {1\over 2}$.
\vskip.5truein
\noindent
\S 2. Proof of Propositions 3 and 1.
\smallskip\noindent
Proposition 1 is an immediate consequence of Proposition 3, since by
a linear change of variable in ${\bf R}^{2n}$, we can assume that
$J(0)=J_{st}$, thus $Q(0)=0$ and $|Q(p)|\leq K |p|^{1\over 2}$.
So we now turn to the Proof of Proposition 3.
\medskip\noindent
2.1. Reduction to Lemma 3. Let $\omega\subset {\bf D}$ be the set of $z\in {\bf D}$ such that
$u\equiv 0$ on a neighborhood of $z$. We need to show that its boundary $b\omega$ is
empty. If it is not empty, we can choose $\zeta_0\in\omega$ such that
${\rm dist}~(\zeta_0,b\omega )=r<1-|\zeta_0|$, and let $\zeta_1\in b\omega$
such that $|\zeta_1-\zeta_0|=r$. So, on the disc defined by $|\zeta-\zeta_0|\leq r$,
$u=0$. We wish to prove that $u\equiv 0$ near $\zeta_1$, getting thus a contradiction.
\smallskip\noindent
Note that the hypotheses of Proposition
3 are preserved under holomorphic change of variable for $z$. Using a conformal map from a neighborhood
of $\zeta_1$ that maps $\zeta_1$ to 0 and that maps the intersection of disc $\{ |\zeta -\zeta_0|<r\}$
with a neighborhood of $\zeta_1$ to a region defined near 0 by $x>-y^2$, and by rescaling, the proof
of Proposition 3 reduces to proving the following:
\vfill\eject
{\bf Lemma 3.}{\it Let $A$ and $\theta >0$ and let $u:{\bf D}\to {\bf C}^n$ be a ${\cal C}^1$ map such that
$$u(x+iy)=0~,~{\rm if}~x\geq -Ay^2~,$$
$$ |{\partial u\over\partial\overline z}|\leq \theta |u|^{1\over 2}|{\partial u\over \partial z}|~.$$
Then $u\equiv 0$ near 0.}
\medskip\noindent
2.2. Further reduction. By replacing $u$ by $Ku$ for $K$ large enough, we can assume that $\theta$
is as small as we wish. In order to apply Lemma 2, we take $\theta={1\over 10}$.
Next (seemingly contradictory to the previous step), we can also assume that
$|u|\leq 1$ and $|{\partial u\over\partial z}|\leq 1$. This can be achieved by rescaling, replacing
the function $z\mapsto u(z)$, by the function $z\mapsto u(\epsilon z))$, for $\epsilon >0$ small
enough. This changes the constant $A$ in the Lemma. We shall therefore use
$\theta ={1\over 10}$, $|u|\leq 1$ and $|{\partial u\over \partial z}|\leq 1$.
\bigskip\noindent
2.2. Proof of Lemma 3. This is the standard game of Carleman's estimates.\hfill\break
In order to apply the $\overline\partial$ estimates of \S 1, we need compact support. Let 
$\chi\in {\cal C}^\infty ({\bf D})$ be such that $0\leq \chi\leq 1$,
$\chi(x+iy)=1$ if $x>-\alpha$, and $\chi (x+iy)=0$ if $x<-2\alpha$, where $\alpha>0$ is 
chosen small
enough so that the region defined by $x<-y^2$ and $x>-2\alpha$ is relatively compact in ${\bf D}$.
\bigskip\noindent
Set $v=\chi u$, so $v$ is a compactly supported map from ${\bf D}$ into ${\bf C}^n$.
We can apply the estimate of Lemma 2, 
with $w(z)={ \partial u\over\partial \overline z}$ if $x\geq -\alpha$ ($z=x+iy$),
and $w(z)=0$ if $x< -\alpha$.
For all $\tau\geq 1$:
$$\int |{\partial \over \partial\overline z}( \chi  u )(z)
-w(z)|^2e^{\tau \phi} \geq
{\tau\over 40} \int |\chi u|^2 e^{\tau \phi}~.$$
However:
\smallskip\noindent
The integrand on the left hand side is zero for $x>-\alpha$. So the left hand side is
at most $O\big( e^{\tau (-\alpha + \alpha^2)\over 2}\big)$. 
If $u(x_0,y)\neq 0$ for some $x_0>\alpha$,
it is immediate to see that (for $\tau$ large)
${\tau\over 80}\int |u|^2 > e^{\tau (-x_0+x_0^2)}$. Letting $\tau\to +\infty$ gives us a contradiction
since $-x_0+{x_0^2\over 2}>-\alpha+{\alpha^2\over 2}$. 
\vskip.5truein
\noindent
\S 3. Examples (proving Proposition 2).
\bigskip\noindent
EXAMPLE: {\it A smooth map $z \mapsto u(z)= (u_1(z),u_2(z))$
from a neighborhood of 0 in ${\bf C}$ into ${\bf C}^2$, such that 
$${|\overline \partial u|\over |d u|}\to 0,~~{\rm as}~~z\to 0~,$$
$u$ vanishes to
infinite order at 0, but is not identically 0 near 0.
\bigskip\noindent
Moreover one can adapt the construction so that $u$ has a non-isolated zero at
0.} 
\vskip.5truein\noindent
3.1. We first give an example without getting a non-isolated zero.
\smallskip\noindent
1) If $n$ is an even positive integer, for $2^{-n}\leq |z|\leq 2^{-n+1}$, set
\smallskip
$$u_1(z)~=~2^{n^2\over 2}z^n$$
$$u_2(z)~=~\chi (z) 2^{(n-1)^2\over 2}z^{n-1}~+~
           (1-\chi (z)) 2^{(n+1)^2\over 2} z^{n+1}~,$$
where $\chi =1$ near $|z|=2^{-n+1}$, $\chi=0$ near $|z|=2^{-n}$,
$|d\chi |=O(2^n)$, and more generally the ${\cal C}^k$ norm of $\chi$
is $O(2^{kn})$.
\bigskip\noindent
2) For $n$ odd, take the same definitions switching $u_1$ and $u_2$.
\bigskip\noindent
CLAIMS:

(a) $u$ is a smooth map vanishing to infinite order at 0, 
\smallskip
(b) For $2^{-n}\leq |z|\leq 2^{-n+1}$
$$|\overline \partial u|\leq {C\over n} |\partial u|~,$$
with $C$ independent on $n$.

\bigskip
\noindent
(a) is straightforward. Note that for $|z|=K 2^{-n}$ (think
${1\over 2}\leq K \leq 2$),\hfill\break
$2^{n^2\over 2}|z^n|=2^{-n^2\over 2}K^n\leq 2^{-n^2\over 3}$.
\bigskip\noindent
We now check (b), for $n$ even. For $n$ odd, the checking
is the same with $u_1$ and $u_2$ switched.
\bigskip\noindent
We have $|du|\geq |du_1|=n2^{n^2\over 2}|z|^{n-1}$.\hfill\break
The factor $n$ will be the needed gain.
\bigskip\noindent
$\overline \partial u_1=0$. So we only have to estimate
$\overline \partial u_2$, in which the non zero terms come
from differentiating $\chi$ whose gradient is of the order of $2^n$.
One gets (with various constants $C$):
$$(\ast )\qquad |\overline \partial u_2|\leq C 
\big(2^n 2^{(n-1)^2\over 2} |z|^{n-1}+2^n 2^{(n+1)^2\over 2} |z|^{n+1}
\big) $$
$$\leq C~|z|^{n-1}\big( 2^{n^2\over 2}+2^{{n^2\over 2}+2n}|z|^2\big) ~.$$
Since $|z|\simeq 2^{-n}$, one indeed gets
$$|\overline\partial u|\leq C 2^{n^2\over 2}|z|^{n-1}~~{\rm so}~~
|\overline \partial u|\leq {C\over n} |du|~~.$$
\bigskip\noindent
{\bf Remark 2.} The example is not difficult but the matter
looks delicate. In particular, the choice of the exponent ${n^2\over 2}$ seems to
be 
somewhat dictated. If in the above definition of $u_1$ we would 
set
$u_1(z)~=~2^{n^2\over p}z^n$ instead of $u_1(z)~=~2^{n^2\over 2}z^n$,
and do the corresponding change in the definition of $u_2$, we
should take $p>1$ for having decay.
Then, for the estimate of the first term on the right 
hand side in ($\ast$), 
we would need $p\leq 2$. But for
the estimate of the second term we would need $p\geq 2$. 
\bigskip
\bigskip\noindent
3.2. We now indicate how to get a non-isolated zero.

\bigskip\noindent
Set $$2^{-n}<r_n={5\over 4}2^{-n}
< a_n={3\over 2}2^{-n}<R_n={7\over 4}2^{-n}<
2^{-n+1}~.$$
We modify the definition of $u_1$ and $u_2$ 
(resp. $u_2$ and $u_1$, if $n$ is odd)
above by setting:
$$u_1(z)~=~2^{n^2\over 2}z^{n-1}(z-a_n)$$
i.e. replacing a factor $z$ by $(z-a_n)$, and accordingly
$$u_2(z)~=~\chi (z) 2^{(n-1)^2\over 2}z^{n-2}(z-a_{n-1})~+~
           \psi (z) 2^{(n+1)^2\over 2} z^n(z-a_{n+1})~,$$
where $\chi =1$ near $|z|=2^{-n+1}$ and $\chi (z)=0$
if $|z|<R_n$ , $\psi=1$ near $|z|=2^{-n}$ and 
$\psi (z)=0$ if $|z|>r_n$, with estimates on the
derivatives, as before. Note that for $r_n<|z|<R_n$, $u_2(z)=0$,
(in particular $u_2(a_n)=0$).
\bigskip\noindent
In the region $r_n<|z|<R_n$, both $u_1$ and $u_2$ are holomorphic.
So the differential inequalities have to be checked only in the regions
$2^{-n}<|z|<r_n$ and $R_n<|z|<2^{-n+1}$. In these regions
$|z|$, $|z-a_n|$, $|z-a_{n+1}|$ and $|z-a_{n-1}|$ all have the same order
of magnitude.
\bigskip\noindent
The estimate for $\overline\partial u_2$ is basically unchanged:
gradient estimates for $\chi$ and $\psi$ and point-wise estimates
of $z^{n-2}(z-a_{n-1})$ and 
$ z^n(z-a_{n+1})$ instead of $z^{n-1}$ and $z^{n+1}$.
\bigskip\noindent
Finally one has to estimate $\partial u_1=2^{n^2\over 2}
\big( (n-1)z^{n-2}(z-a_n)+z^{n-1}\big)$. For $n$ large, in the regions under consideration
the first term in the parenthesis dominates the second one
($|z-a_n|\geq {1\over 8}|z|$),
and one has $|\partial u_1|\geq {n\over 10}2^{n^2\over 2}|z|^{n-1}$.
\bigskip\noindent
So, as previously 
$$|\overline \partial u|\leq {C\over n} |\partial u|~,$$
with $C$ independent on $n$.\hfill\break
We have $u_1(a_n)=u_2(a_n)=0$.
\bigskip\noindent
{\bf Remark 3.} Almost complex structures that are merely continuous and not
H\"older continuous are certainly of much less interest. However for these
merely continuous almost complex structures, on can still
prove the existence of many
$J$-holomorphic curves, see 5.1 in [9]. The above example, 
by adding four real dimensions,
allows one to find a continuous almost complex structure $J$ on ${\bf C}^4$,
such that there is a smooth $J$-holomorphic map $U:{\bf D}\to ({\bf C}^4,J)$
non identically 0 near 0, but vanishing to infinite order at 0, 
and with a non-isolated zero. We are still very far from an example with
a H\"older continuous almost complex structure, as asked in [10].
\smallskip\noindent 
We set $U(z)=(u(z),u_3(z),u_4(z))$,
where $u=(u_1,u_2)$ is the map of the above example, with $u(a_n)=0$. 
Looking at the proof, we see that
$\partial u(z) \neq 0$ unless $z=0$ or, 
$z={n-1\over n}a_n$, for some $n$, and that $\overline \partial u=0$ near ${n-1\over n}a_n$. 
For each $n$, let $\psi_n$ be a non negative smooth function with support in a 
small neighborhood of $\{ 2^{-n}\leq |z|\leq 2^{-n+1}\}$, such that: 
$\psi_n(z)=|z-a_n|^2$ near $a_n$, $\psi_n(z)>0$ if $2^{-n}\leq |z|\leq 2^{-n+1}$ and $z\neq
a_n$, $\psi_n$ is constant on a neighborhood of ${n-1\over n}a_n$. Note that
$\overline\partial \psi_n(a_n)=0$, and $\overline\partial \psi_n = 0$ near ${n-1\over n}a_n$.
For $\epsilon_n$ small enough, set $u_3=\sum_n\epsilon_n\psi_n$.
Then $u_3$ is a smooth 
function that
vanishes only at the points $a_n$ and at $0$, such that $\overline\partial u_3=0$
at $a_n$ and near ${n-1\over n}a_n$, and such that   
${|u_3|+|\nabla u_3|\over |\partial u|}\to 0$, as $z\to 0$. Finally one sets
$u_4(z)=zu_3(z)$. Note the following injectivity property: $(u_3(z),u_4(z))=
(u_3(z'),u_4(z'))$ if and only if $z=z'$, unless $z=a_n$ for some $n$ or $z=0$.
One can define the almost
complex structure on ${\bf C}^4$
by defining the matrix $Q$ in $(E)$. We set $Q(0)=0$ and we
need
$[Q(u_1(z),u_2(z),u_3(z), u_4(z))](\partial U)
=\overline{\overline\partial U}$. 
Due to the injectivity
property of $U$ and the vanishing of $\overline \partial U$ near the
points ${n-1\over n}a_n$ where $\partial u = 0$, and at the points $a_n$
where $U=0$, the above requirement on $Q$ 
is compatible with the 
requirement of continuity and $Q(0)=0$, since ${|\overline\partial U|\over |\partial u|}(z)\to
0$ as $z\to 0$. 
\smallskip
\noindent
\S 4. Appendix.
\smallskip
\noindent
4.1. The inequality (IN1) has been used by many authors 
([3] and [14] Lemma 3.2.4 - see also [4], [8], [12], [13]). 
Given a bounded matrix valued function $z\mapsto A(z)$, 
where $A(z)$ is a $n\times n$ matrix,
one can solve the equation
${\partial M\over\partial \overline z}+AM$ with solution $M(z)$ an invertible matrix.
Locally this is easily obtained by the inverse function Theorem, for a global result
see [11] (or the Appendix p. 61 in [13]). 
If (IN1) is satisfied, there exists a bounded matrix valued map $z\mapsto A(z)$
(depending on $u$) such that
${\partial u\over\partial\overline z}+A(z)u=0$. 
Define $v$ by $u=Mv$.
The above equation yields ${\partial v\over \partial \overline z}=0$,
so $v$ is holomorphic, and the zero set of $u$ is the same as the zero set of the
holomorphic vector valued function $v$. 
Globally on ${\bf D}$, that will give the Blaschke
condition for the zero set of $J$-holomorphic maps, if $J$ is ${\cal C}^1$, and
if $|\nabla u|$ is bounded.
See more on that topic in [9].  
\bigskip\noindent
4.2. The scalar case of (IN2) is $|{\partial u\over \partial \overline z}|\leq
a |{\partial u\over \partial z}|$ with $u$ scalar valued, and we should take $0\leq a<1$. 
For each function $u$, satisfying the inequality, there is a
bounded measurable (not continuous) function $\alpha$, with $|\alpha (z)|\leq a<1$, such that
$u$ satisfies the Beltrami equation  
${\partial u\over \partial \overline z}=\alpha (z) {\partial u\over \partial z}$.
Beltrami equations have been
much studied (Bers, Bers-Nirenberg, Morrey, Vekua, $\cdots$). A very convenient 
reference is [2]. It is shown in [2], Theorems 3.1 and 3.2, that,
given $\alpha$ (in our case, not given a priori but associated to $u$),
with $|\alpha |\leq a<1$,  
there exists a H\"older continuous change of variables $\psi$ such that 
a function $v$
satisfies the Beltrami equation ${\partial v\over\partial \overline z}=\alpha(z)
{\partial v\over\partial z}$
if and only if $v\circ \psi$ is holomorphic. So there is again
reduction to the holomorphic case, by composition on the right rather than
by composition on the left as in 4.1. However, this reduction is only in the scalar case and it
is much more difficult. Note that, even if we started from the point of view of
$J$-holomorphicity ${\partial u\over\partial \overline z}=
\overline{\beta (u) {\partial u\over \partial z}}$, we switched here to the
quasi-conformal point of view ${\partial v\over\partial\overline z}=
\alpha (z) {\partial v\over \partial z}$. This, in the theory of almost complex
structures on ${\bf C}$, corresponds to
studying maps from $({\bf C},J_{st})$ to $({\bf C},J)$, or vice versa.
And it corresponds to the non-linear and linear approaches to the Theorem of
Newlander-Nirenberg.

\smallskip

\centerline{REFERENCES.}
\smallskip\noindent
\item{[1]} S.\ Alinhac. Non-unicit\'e pour des op\'erateurs
diff\'erentiels \`a caract\'eristiques complexes simples.
Ann. Sci. \'Ecole Norm. Sup. {\bf 13} (1980), no. 3, 385-393.

\item{[2]} B.V.\ Bojarski. Generalized solutions of a system
of differential equations of the first order and elliptic type
with discontinuous coefficients. University of Jyv\"askyl\"a.
Department of Mathematics and Statistics Report 118 (2009).
(Translation of a Math. Sbornik paper of 1957).

\smallskip
\item{[3]} A.\ Floer, H.\ Hofer, and D.A.\ Salamon. Transversality 
in elliptic Morse theory for the symplectic action. Duke Math. J. 
{\bf 80} (1995), no. 1, 251-292.

\smallskip    
\item{[4]} X.\ Gong, J.P.\ Rosay. Differential inequalities of
continuous functions and removing singularities of Rado type for
$J$-holomorphic maps. Math. Scand. {\bf 101} (2007), 293-319.

\smallskip
\item{[5]} L.\ H\"ormander. {\it The Analysis of Linear Partial
Differential Operators I.} Grundlehren der Math. Wissenschaften
{\bf 256} (1983) Springer-Verlag.

\smallskip
\item{[6]} L.\ H\"ormander. {\it The Analysis of Linear Partial
Differential Operators III.} Grundlehren der Math. Wissenschaften
{\bf 274} (1985) Springer-Verlag.

\smallskip
\item{[7]} L.\ H\"ormander. {\it Notions of Convexity.} Progress in Math.
{\bf 127} (1994) Birkh\"auser.

\smallskip
\item{[8]} S.\ Ivashkovich, J.P.\ Rosay.  Schwarz-type lemmas for solutions
of $\overline \partial$-inequalities and complete hyperbolicity of almost
complex manifolds. Ann. Inst. Fourier 54 (2004), no. 7, 2387-2435.

\smallskip
\item{[9]} S.\ Ivashkovich, J.P.\ Rosay. Boundary values and boundary
uniqueness of \hfill\break $J$-holomorphic mappings. 
arXiv:0902.4800. 

\smallskip
\item{[10]} S.\ Ivashkovich, V.\ Shevchishin. Local properties of
J-complex curves in Lipschitz-continuous structures. arXiv:0707.0771.

\item{[11]} B.\ Malgrange. {\it Lectures on the theory of several 
complex variables. Notes by Raghavan Narasimhan}. Tata Institute of
Fundamental Research, Bombay (1958).

\smallskip
\item{[12]}  D.\ McDuff - D.\ Salamon. {\it $J$-holomorphic curves and 
quantum cohomology}. Univ.
Lect. Series AMS, 6 (1994).

\smallskip
\item{[13]} J.P.\ Rosay. Notes on the Diederich-Sukhov-Tumanov normalization 
for almost complex structures. Collect. Math. {\bf 60} (2009), no. 1, 43-62.

\smallskip
\item{[14]} J.C.\ Sikorav. Some properties of holomorphic curves in
almost complex manifolds. In {\it Holomorphic Curves in Symplectic
Geometry}, eds. M. Audin and J. Lafontaine, Birkhauser (1994), 351--361.
\bigskip
\noindent
Jean-Pierre Rosay

\noindent
Department of Mathematics,
University of Wisconsin, Madison WI 53705.

\noindent
jrosay@math.wisc.edu

\bigskip\noindent
2000-AMS classification: Primary 32Q65. Secondary 35R45, 35A02.\hfill\break
Key words: J-holomorphic curves, differential inequalities, uniqueness.

\end